\documentclass[oneeqnum,onefignum,onetabnum,onealgnum]{siamart251216}
\usepackage{amssymb,mathtools}
\usepackage{microtype}
\AddToHook{cmd/normalsize/after}{%
  \setlength{\abovedisplayskip}{6pt plus 1pt minus 2pt}%
  \setlength{\belowdisplayskip}{6pt plus 1pt minus 2pt}%
  \setlength{\abovedisplayshortskip}{0pt plus 1pt}%
  \setlength{\belowdisplayshortskip}{4pt plus 1pt minus 1pt}%
}
\AtBeginDocument{\normalsize\setlength{\jot}{2pt}\raggedbottom}

\usepackage{booktabs,graphicx,xcolor,longtable,pdflscape}
\usepackage{algpseudocode}
\usepackage{enumitem,needspace,float}
\hypersetup{hypertexnames=false,colorlinks=true,linkcolor=blue,citecolor=blue,urlcolor=blue}
\newsiamthm{example}{Example}
\newsiamremark{remark}{Remark}
\theoremstyle{plain}
\theoremheaderfont{\normalfont\scshape}
\theorembodyfont{\normalfont}
\theoremseparator{.}
\theoremsymbol{}
\newtheorem{assumption}{Assumption}
\headers{Quasi-Newton Methods on Matrix Manifolds}{Z. Tian, J. Hu, Z. Yao, and B. Li}
\newcommand{\rev}[1]{#1}
\newcommand{\revdisplay}[1]{#1}
\newcommand{\companionrev}[1]{#1}
\newcommand{\kyrev}[1]{#1}
\newcommand{\R}{\mathbb{R}}\newcommand{\M}{\mathcal{M}}\newcommand{\Om}{\Omega}
\DeclareMathOperator{\grad}{grad}\DeclareMathOperator{\Hess}{Hess}\DeclareMathOperator{\dist}{dist}
\DeclareMathOperator{\Skew}{Skew}\DeclareMathOperator{\St}{St}\DeclareMathOperator{\OB}{OB}
\DeclareMathOperator{\SO}{SO}

\newcommand{\norm}[1]{\lVert #1\rVert}
\title{A Euclidean Quasi-Newton Method for Matrix-Manifold Optimization via Almost-Global Coordinates\thanks{An earlier version of this manuscript is available as arXiv:2611.24101.}}
\author{Ziang Tian\thanks{Qiuzhen College, Tsinghua University
(\email{tza23@mails.tsinghua.edu.cn}).}
\and Jiang Hu\thanks{Yau Mathematical Sciences Center, Tsinghua University
(\email{jianghu@tsinghua.edu.cn}).}
\and Zhigang Yao\thanks{Corresponding author. Department of Statistics and Data Science,
National University of Singapore (\email{zhigang.yao@nus.edu.sg}).}
\and Bingjie Li\thanks{Corresponding author. School of Public Health, Tongji Medical College,
Huazhong University of Science and Technology (\email{bingjieli@hust.edu.cn}).}}
\hypersetup{pdfauthor={Ziang Tian, Jiang Hu, Zhigang Yao, Bingjie Li}}
\begin{document}\maketitle
\begin{abstract}
We study smooth matrix-manifold optimization using objective values and
first derivatives. \rev{Almost-global coordinates transform} the
manifold-constrained problem into an unconstrained pullback problem, to
which we apply BFGS with a chart-restart strategy and a limited-memory
implementation. We establish correspondence of first- and second-order
stationary points and minimizers, together with conditional convergence
and local $q$-superlinear convergence. The local rate follows
from Euclidean full-memory BFGS theory once the chart is fixed, without
the locking condition required by classical transported-gradient
Riemannian BFGS analyses. Tests on 15 benchmark problems from Manopt and
Pymanopt, with ten instances per problem, show that \rev{our} limited-memory
method solves all 150 cases, \rev{whereas} Riemannian BFGS and Riemannian
steepest descent \rev{solve 132 and 117 cases, respectively}. On the 115 cases
\rev{successfully} solved by all three methods, \rev{our method takes roughly half
the time of either baseline, measured by} geometric-mean solve time.

\end{abstract}

\section{Introduction}\label{sec:intro}
We consider the smooth optimization problem
\begin{equation}\label{eq:P}
\min_{x\in\M} f(x),\tag{P}
\end{equation}
where $\M$ is a matrix manifold---a smooth submanifold of a matrix space, or a quotient thereof---and
$f:\M\to\R$ is twice continuously differentiable. Such problems have
important applications in numerical linear algebra, statistical estimation,
signal processing, and machine learning~\cite{AMS2008,Boumal2023,HuLiuWenYuan2020}.

Different applications give rise to different matrix manifolds.
In eigenvalue computation and subspace estimation, orthonormal bases
are represented on the Stiefel manifold, while the corresponding
subspaces are points on the Grassmann manifold~\cite{EAS1998}.
Their objectives often take a trace form, such as
$f(X)=-\operatorname{tr}(X^\top AX)$ for symmetric $A$.
Related subspace problems arise in manifold fitting~\cite{YaoXia2025},
where local tangent or normal subspaces describe the geometry to be estimated.
With a fixed center and neighborhood weights, smooth robust subspace-fitting
losses give objectives on the Grassmann manifold to which the proposed
coordinate-based quasi-Newton method applies\companionrev{; see the
illustrative local fitting objective in~\cite[Section 7.1]{GP}}.
In shape alignment and orthogonal Procrustes analysis, the unknown
transformation belongs to the orthogonal group~\cite{Schonemann1966}.
Covariance estimation and matrix averaging lead to the manifold of
symmetric positive-definite matrices~\cite{SraHosseini2015}.
Low-rank matrix completion and sensing exploit fixed-rank manifolds
to represent the unknown matrix~\cite{Vandereycken2013}.
Frame design and semidefinite relaxations provide further examples
involving products of spheres and oblique manifolds~\cite{HuLiuWenYuan2020}.

Orthogonal and unitary constraints are also used in recurrent neural
networks to control the amplification of hidden states and gradients;
the objective is then the training loss of the network~\cite{Arjovsky2016,Wisdom2016}.
Cayley and exponential parametrizations provide concrete ways to enforce
these constraints~\cite{Helfrich2018,LiFuxin2019,Lezcano2019}.
These applications motivate algorithms that exploit first-order information
while limiting the cost of geometric operations.

\subsection{Literature Review}\label{ssec:riemannian}
The relevant algorithmic class is defined by its derivative requirements:
methods may evaluate $f$ and its first derivative, but do not require an
objective Hessian or a Hessian--vector oracle. Efficiency involves the
number and cost of these oracle calls, the cost of geometric operations,
and memory use. Gradient descent, nonlinear conjugate gradients, and
quasi-Newton methods belong to this class~\cite{NW2006,SunYuan2006}. BFGS builds a curvature
approximation from gradient differences, so it is a quasi-Newton method
using first-order information, even though its local convergence can be
superlinear. \rev{Exact-Hessian Newton and trust-region methods can be
expensive or infeasible when second derivatives are costly to evaluate.}

\rev{The geometric foundation is the retraction-based
framework~\cite{AMS2008}, developed further
in~\cite{Boumal2023,SatoBook}.} Riemannian gradient
methods move along retractions, while conjugate-gradient methods also
transfer previous directions between tangent spaces. This framework
makes both first-order oracle work and the cost of geometric primitives
relevant when comparing algorithms.
\rev{Spectral gradient methods use } Barzilai--Borwein
\rev{stepsizes}~\cite{BarzilaiBorwein1988,Raydan1997}. \rev{Related step-selection
strategies have also been developed} for Riemannian optimization and
orthogonality-constrained problems~\cite{IannazzoPorcelli2018,WenYin2013}.

Nonlinear conjugate gradients (CG) use a scalar recurrence to retain
directional information, with low storage cost. In the Riemannian setting,
a typical update is
\[
 d_k=-g_k+\beta_k^{\mathrm{CG}}\mathcal T_{k-1}d_{k-1},
 \qquad g_k=\grad f(x_k),
\]
where $\mathcal T_{k-1}$ moves the previous direction to $T_{x_k}\M$.
The coefficient and line-search conditions must be designed together;
\rev{The Euclidean method in~\cite{DaiYuan1999} provides an important
globally convergent construction under Wolfe conditions.
This direction has also been developed on manifolds.
Scaled vector transport is introduced in~\cite{SatoIwai2015} to
relax the nonexpansiveness restriction used in global convergence proofs.
A Dai--Yuan-type Riemannian CG method using weak Wolfe conditions
is established in~\cite{Sato2016DY}. The later framework
in~\cite{Sato2022} unifies several transport- and
inverse-retraction-based CG constructions and provides convergence
analyses for specific choices.} These results
show that transport scaling and line-search compatibility can yield
effective first-order methods without requiring the full quasi-Newton
transport structure. Their global convergence statements concern
stationarity under assumptions, rather than a blanket superlinear rate.

Quasi-Newton methods seek faster local convergence using the same
first-order oracle. Classical full-memory BFGS accumulates secant
information to approximate curvature~\cite{NW2006}; L-BFGS retains only
$\ell$ pairs and computes its direction by a two-loop recursion, with
$O(\ell d)$ storage and arithmetic for a $d$-dimensional coordinate
problem, excluding oracle and chart costs~\cite{LiuNocedal1989}.
The target here is chart-restarted Euclidean BFGS, with an L-BFGS
implementation for the experiments. In a fixed chart, its direction is
$p_k=-H_k\nabla h(y_k)$ for $h=f\circ G$, and its secant pair is
\[
 s_k=y_{k+1}-y_k,\qquad
 z_k=\nabla h(y_{k+1})-\nabla h(y_k).
\]
Both vectors live in the same Euclidean space. The implementation combines
Wolfe line search, cautious updates, Riemannian gradient stopping criteria, and
radius-based recentering (\S\ref{sec:bfgs-final-design}).

To state the local rate precisely, $q$-superlinear convergence means
\[
 \lim_{k\to\infty}
 \frac{\dist(x_{k+1},x^\star)}{\dist(x_k,x^\star)}=0.
\]
 It is an asymptotic statement,
not a bound on total runtime or a consequence of a small final residual.
For full-memory BFGS, standard smoothness, positive-definite local Hessian,
and line-search hypotheses support a Dennis--Mor\'e analysis~\cite{NW2006,DennisMore1974}.
Our local \rev{rate applies once} chart changes stop (\rev{Lemma}~\ref{thm:super})\rev{; conditions under which chart restarts stop are given in Theorem}~\ref{prop:restart-finite}. Exact within-chart pairs remove the transport
issue, but do not alone prove this rate. Ordinary fixed-memory L-BFGS does
not automatically inherit the full-memory theorem; methods that preserve
full-BFGS information by displacement aggregation provide a distinct
route to limited-storage superlinear guarantees~\cite{BerahasCurtisZhou2022}.
The fixed-memory experiments in this paper are therefore reported as
performance evidence, separately from the full-memory local theorem.

\rev{For geometric BFGS, convergence is analyzed
in~\cite{RingWirth2012} using compatibility conditions between
retractions and transports.
In the transported-gradient Broyden framework
of~\cite{HuangGallivanAbsil2015}, isometry and locking are two distinct
requirements.} Isometry preserves tangent-space inner products:
\[
 g_{R_x(\eta)}(\mathcal T_\eta u,\mathcal T_\eta v)=g_x(u,v).
\]
Locking aligns the transported step with the line-search velocity:
\[
 \mathcal T_\eta\eta=\beta_\eta\,DR_x(\eta)[\eta],\qquad
 \beta_\eta=\frac{\norm{\eta}_x}{\norm{DR_x(\eta)[\eta]}_{R_x(\eta)}}.
\]
For that framework's scaled gradient-difference pair, these conditions
connect Wolfe curvature to positive secant curvature and maintain positive
definiteness. They restrict the permissible retraction--transport pair;
preserving norms alone does not establish locking.

These are strict requirements of that classical formulation, rather than
universal necessary conditions for every Riemannian BFGS variant.
\rev{The method in~\cite{HuangAbsilGallivan2018} uses cautious updates
and a weak Armijo--Goldstein line search to avoid differentiated-retraction
information and the associated locking requirement, while establishing
conditional local superlinear convergence.} Transport compatibility still
matters in the local analysis. Accordingly, failure to satisfy classical
locking does not by itself prove loss of superlinear convergence or
numerical stagnation. The question motivating our formulation is how to
obtain effective curvature learning from first derivatives with less
transport construction and validation.

\subsection{Contribution}\label{ssec:contrib}
Our contributions are twofold.
\begin{itemize}[leftmargin=1.5em,itemsep=2pt]
\item \kyrev{We develop chart-restarted Euclidean BFGS and a limited-memory
implementation for matrix manifolds using almost-global coordinates~\cite{GP}.
We quantify how the chart radius controls Jacobian conditioning for several
geometries, including Stiefel and Grassmann, and how the Jacobian condition
number enters the pullback Hessian bound at a nondegenerate minimizer
(\S\ref{sec:chart-restart}).}
On 15 benchmark problems from Manopt and Pymanopt with ten instances
each, the limited-memory method solves all 150 cases at a Riemannian gradient
tolerance $10^{-6}$, compared with 132 for Riemannian BFGS and 117 for
steepest descent. On the 115 cases solved by all three methods, its
geometric-mean solve time is \rev{roughly half that of either baseline}
(\S\ref{sec:num}).
\item \kyrev{We establish optimality correspondence within the chart image
and equality of infima (Theorem~\ref{thm:corr}). Under the assumptions of
Theorem~\ref{prop:restart-finite}, a chart centered sufficiently near a
nondegenerate minimizer needs no further radius restart, and full-memory
BFGS terminates at the minimizer or converges locally $q$-superlinearly.
The rate follows from fixed-coordinate Euclidean BFGS theory without
vector transport or a locking condition.}
\end{itemize}

Section~\ref{sec:manifolds} introduces the pullback formulation and the
\rev{almost-global coordinate} preliminaries, then gives the maps for the test
geometries. Section~\ref{sec:alg} develops the algorithm and convergence
analysis. Section~\ref{sec:num} reports the experiments, and
Section~\ref{sec:conc} concludes.

\section{Preliminaries}\label{sec:manifolds}
We use the Frobenius inner product on ambient matrix spaces.
The Riemannian gradient and Hessian are denoted by $\grad f$ and
$\Hess f$; hats denote their representations in orthonormal tangent
frames, as defined in Section~\ref{sec:pullback-derivatives}.
Riemannian gradient norms use the specified metric, while coordinate
gradients use the Euclidean norm. For an ambient extension $\bar f$,
we write $\nabla\bar f$.

\subsection{Reformulation in coordinates}\label{ssec:prior}
Following early pullback and matrix-parametrization
approaches~\cite{Dreisigmeyer2007,Arsigny2007,HylandRatsch2017}, a feasible
map $G:\Om\subseteq\R^d\to\M$ yields the unconstrained pullback problem
\begin{equation}\label{eq:H}
 \min_{y\in\Om}h(y),\qquad h=f\circ G.\tag{H}
\end{equation}
BFGS uses function values and gradients of $h$, without evaluating an
objective Hessian. We use the almost-global coordinate framework
of~\cite{GP}\companionrev{, using its strict approximate parametrizations}.
\companionrev{The coordinate map and its inverse are exact; the omitted
set is the only sense in which the parametrization is approximate.
We equip the parameter domain with the Euclidean metric, rather than
the pullback metric used for the gradient identity in~\cite[Proposition 2.7]{GP}.}
Earlier moving-base-point and dynamic-trivialization
methods~\cite[Section 6]{Dreisigmeyer2018}\cite{Lezcano2019} motivate
recentring when coordinates become unsuitable; scaled Cayley and
exponential representations provide related orthogonality-preserving
constructions~\cite{Helfrich2018,LezcanoMR2019}.
Our method accumulates quasi-Newton information in fixed coordinates
and clears it on restart, avoiding per-step tangent-vector transport.
The superlinear analysis concerns full-memory BFGS; the experiments use
a limited-memory implementation.

\subsection{Almost-global coordinates}
\companionrev{We use the $C^2$ version of the strict approximate
parametrization in~\cite[Definition 2.2]{GP}.}
\begin{definition}[Almost-global coordinates~\companionrev{\cite[Definition 2.2]{GP}}]\label{def:am}
\companionrev{A $C^2$ Riemannian $d$-manifold $\M$ admits almost-global
coordinates if there is a $C^2$ diffeomorphism
$G:\Om\to\M\setminus\Sigma$, where $\Om\subseteq\R^d$ is open and
$\Sigma$ is closed and has Riemannian volume zero. The parameter domain
may be disconnected. The map $G$ is the coordinate representation,
and $y=G^{-1}(x)$ gives coordinates on its open dense image.}
\end{definition}
\begin{lemma}[Products and diffeomorphisms~\companionrev{\cite[Theorem 3.1 and Section 2.1]{GP}}]\label{lem:prod}
Finite products and $C^2$ diffeomorphic images preserve almost-global
coordinates. For product charts the omitted set is
$\bigcup_i\pi_i^{-1}(\Sigma_i)$; under a diffeomorphism $\Phi$ it is
$\Phi(\Sigma)$.
\end{lemma}
These properties apply componentwise when a manifold is disconnected.

\subsection{Coordinates for the test manifolds}\label{ssec:test-charts}
The Stiefel construction contains the sphere and rotation-group cases.
We recall its almost-global property from~\cite{GP}, then describe the
product extensions and the graph coordinates used for Grassmann.
\companionrev{These maps and their inverses are verified directly in
\cite[Section 6]{GP}; a null omitted set in a group does not, by itself,
give a null omitted set in a quotient.}
Write $\Skew(n)$ for skew-symmetric matrices, $\St(n,p)$ for the Stiefel
manifold, and $\SO(n)$ for the special orthogonal group.

\begin{proposition}[Stiefel coordinates~\companionrev{\cite[Theorem 6.1 and Corollary 6.3]{GP}}]\label{prop:stiefel}
Let $1\le p\le n$, $X_0\in\St(n,p)$, and
$Q_0=[X_0\;X_0^\perp]\in\mathrm O(n)$.
For $A\in\Skew(p)$ and $B\in\R^{(n-p)\times p}$, set
$M(A,B)=I+A+B^\top B$. The map
\begin{equation}\label{eq:stiefel-chart}
 G_{\St}(A,B)=Q_0\begin{bmatrix}2M(A,B)^{-1}-I\\
 2B M(A,B)^{-1}\end{bmatrix}
\end{equation}
is a diffeomorphism onto $\{X\in\St(n,p):\det(I+X_0^\top X)\ne0\}$
when $p<n$, with closed, volume-null complement.
For $p=n$ it gives almost-global coordinates on $X_0\SO(n)$;
the two components of $\mathrm O(n)$ require separate charts.
The case $p=1$ recovers inverse stereographic coordinates on
$S^{n-1}\setminus\{-X_0\}$.
\end{proposition}

\kyrev{For $p<n$, the real Stiefel map in~\eqref{eq:stiefel-chart} is the inverse generalized
left-localized Cayley map of Kume and Yamada~\cite[Proposition 2.2]{KumeYamadaCayley},
with their center $S=Q_0$ and blocks $V_{11}=A$, $V_{21}=-B$.
Their structured centers also yield efficient formulas for the map, its inverse,
and the pullback gradient~\cite[Sections 2.2--2.3]{KumeYamadaCayley}.}

By Lemma~\ref{lem:prod}, columnwise sphere charts give coordinates on
$\OB(d,m)=(S^{d-1})^m$.
Whitening gives $G_B=B^{-1/2}G_{\St}$ on
$\St_B(n,p)=\{X:X^\top BX=I\}$ for $B\succ0$.
Products with identity coordinates give the factor representation
$X=UCV^\top$ on $\St(m,r)\times\R^{r\times r}\times\St(n,r)$.
\companionrev{Complex analogues use adjoints and real Frobenius
products~\cite[Proposition 6.4]{GP}.}

\begin{proposition}[Grassmann graph coordinates~\companionrev{\cite[Theorem 6.5]{GP}}]\label{prop:grassmann}
Let $1\le p<n$ and $Q_0=[X_0\;X_0^\perp]\in\mathrm O(n)$.
The map
\begin{equation}\label{eq:grassmann-chart}
 G_{Q_0}(B)=\operatorname{colspan}\!\left(
 Q_0\begin{bmatrix}I_p\\B\end{bmatrix}\right),
 \qquad B\in\R^{(n-p)\times p},
\end{equation}
is a diffeomorphism onto the open dense set of subspaces with invertible
projection onto $\operatorname{colspan}(X_0)$. Its complement is closed
and has Riemannian volume zero.
\end{proposition}

Matrices $X$ and $XO$, $O\in\mathrm O(p)$, represent the same Grassmann
point. The implementation uses thin QR of $[I_p;B]$, followed by
column-sign correction and multiplication by $Q_0$, to select a representation matrix.

\section{Algorithm and convergence analysis}\label{sec:alg}
\rev{We first derive the pullback derivatives and establish the
correspondence of optimality conditions. We then present BFGS in a fixed
chart and the chart-restart strategy. Finally, we analyze convergence
at restart centers and local superlinear convergence once the chart
is fixed.}

\subsection{The pullback and its derivatives}\label{sec:pullback-derivatives}
\rev{Let} $f\in C^2(\M)$, and \rev{let} $G$ \rev{be an almost-global coordinate
representation of} $\M$ (Definition~\ref{def:am})\rev{. The optimality results
below use these smoothness and invertibility properties; compactness
conditions for convergence will be stated separately}. At $x=G(y)$,
choose a $g_x$-orthonormal frame $E_x:\R^d\to T_x\M$, and let
$E_x^*:T_x\M\to\R^d$ be its metric adjoint. We use
\[
 \widehat{\grad f}(x):=E_x^*\grad f(x),\qquad
 \widehat{\Hess f}(x):=E_x^*\circ\Hess f(x)\circ E_x,
 \qquad J_y:=E_x^*\circ dG_y.
\]
Thus the hatted quantities are tangent-frame coordinates, whereas the
unhatted quantities are tangent vectors or operators on the manifold.
In particular, $J_y$ is an invertible $d\times d$ matrix. The maps and metrics in Section~\ref{sec:manifolds} define this
differential explicitly. In computation, differentiating $G$ and
expressing its tangent output in an orthonormal frame gives $J_y$;
its singular values and condition number are numerically computable.
Equivalently, the eigenvalues of the pullback metric matrix
$J_y^\top J_y$ give its squared singular values. The next lemma relates
the Euclidean gradient and Hessian to the Riemannian gradient, Hessian,
and metric through the chart differential.
\begin{lemma}\label{lem:deriv}
Let $h=f\circ G$ and $x=G(y)$. Then
$\nabla h(y)=J_y^\top\widehat{\grad f}(x)$. For $u,v\in\R^d$,
\[
 D^2h(y)[u,v]
 =\Hess f(x)[dG_yu,dG_yv]
 +g_x\bigl(\grad f(x),(\nabla dG)_y[u,v]\bigr),
\]
where $\nabla dG$ is the covariant derivative of the differential, using
the Euclidean connection on $\Om$ and the Levi--Civita connection on $\M$.
At a stationary point this reduces to
$\nabla^2h(y)=J_y^\top\widehat{\Hess f}(x)J_y$.
\end{lemma}
The following theorem establishes the correspondence between stationary
points and local minimizers of the original problem \eqref{eq:P} and the
unconstrained problem \eqref{eq:H} \rev{within} the \rev{chart} image.
\begin{theorem}[Correspondence]\label{thm:corr}
$G$ is a bijection between the stationary points of $h$ and the critical points of $f$ in $\M\setminus\Sigma$;
$\nabla h(y)=0\iff\grad f(G(y))=0$. At such points $\nabla^2 h(y)=J_y^\top\widehat{\Hess f}(x)J_y$, so the
correspondence preserves the Morse index and nondegeneracy. In particular,
the first-order stationary points correspond, and
\[
 \nabla h(y)=0,\quad \nabla^2h(y)\succeq0
 \quad\Longleftrightarrow\quad
 \grad f(x)=0,\quad \Hess f(x)\succeq0.
\]
Thus second-order stationary points correspond as well.
Furthermore, $y$ is a (strict) local minimizer of $h$ if and only if
$G(y)$ is a (strict) local minimizer of $f$. Moreover $\inf_\Om h=\inf_\M f$, and if a global
minimizer of $f$ lies in $\M\setminus\Sigma$ it corresponds to a global minimizer of $h$.
\end{theorem}
\begin{proof}
By Lemma~\ref{lem:deriv} and invertibility of $J_y$, the gradients vanish together; $G$ being a bijection onto
$\M\setminus\Sigma$ transfers this to the stationary sets. At a stationary point the gradient term drops and
$\nabla^2h=J_y^\top\widehat{\Hess f}J_y$ is a congruence, preserving inertia. Local minimality, including strict local minimality, is preserved because
$G$ is a diffeomorphism onto the open set $\M\setminus\Sigma$.
Finally $G$ is a bijection
$\Om\to\M\setminus\Sigma$, so $\inf_\Om h=\inf_{\M\setminus\Sigma}f=\inf_\M f$ by density of $\M\setminus\Sigma$
and continuity of $f$.
\end{proof}

\subsection{Chart Restart BFGS}\label{sec:bfgs-final-design}
We first describe BFGS in fixed coordinates and its numerical
implementation, then develop a chart-restart framework motivated by
coordinate conditioning.

\subsubsection{BFGS in a fixed chart}\label{sec:fixed-chart-bfgs}
By Theorem~\ref{thm:corr}, the pullback $h=f\circ G_x$ preserves the
relevant optimality properties. Write $g_j=\nabla h(y_j)$ and choose a
Wolfe step $y_{j+1}=y_j-\alpha_jH_jg_j$, with $H_j\succ0$.
Algorithm~\ref{alg:fixed-bfgs} gives the standard BFGS update.
\kyrev{Within a fixed chart, the secant pairs lie in one Euclidean space,
as discussed for Cayley coordinates by Kume and Yamada~\cite[Remark 3.6]{KumeYamadaCayley}.}
\begin{algorithm}[tbp]
\caption{BFGS in a fixed chart}\label{alg:fixed-bfgs}
\begin{algorithmic}[1]
\Require A chart $G_x$, $y_0=0$, $H_0\succ0$, and tolerance $\tau>0$
\State Set $h=f\circ G_x$ and $g_0=\nabla h(y_0)$
\For{$j=0,1,\ldots$}
\If{$\norm{\grad f(G_x(y_j))}<\tau$}
\State \Return $G_x(y_j)$
\EndIf
\State Choose a Wolfe step $\alpha_j$ along $p_j=-H_jg_j$
\State $y_{j+1}=y_j+\alpha_jp_j$; $g_{j+1}=\nabla h(y_{j+1})$
\State $s_j=y_{j+1}-y_j$; $z_j=g_{j+1}-g_j$
\State $\theta_j=(s_j^\top z_j)^{-1}$
\State $H_{j+1}=(I-\theta_js_jz_j^\top)H_j(I-\theta_jz_js_j^\top)+\theta_js_js_j^\top$
\EndFor
\end{algorithmic}
\end{algorithm}

\rev{Since $h(y)=f(G_x(y))$ and the two problems have the same infimum,
solving $h$ to a prescribed objective accuracy yields the same objective
accuracy for $f$ at $G_x(y)$.}
\companionrev{Equality of infima follows from density and continuity
\cite[Theorem 2.8]{GP}; it does not ensure attainment in the chart,
since a manifold minimizer may belong to the omitted set.}

\subsubsection{Conditioning-driven chart restart}\label{sec:chart-restart}
Optimality is preserved by a chart, but conditioning need not be.
Uneven singular values of its differential can slow the Euclidean
iteration. Dynamic recentering~\cite{Lezcano2019,LezcanoMR2019} addresses
this by rebuilding $G_{x_{k+1}}$ at $x_{k+1}=G_{x_k}(y)$, clearing
secant history, and initializing the next solve at zero with a positive
scalar inverse-Hessian scale.
\kyrev{Kume and Yamada~\cite{KumeYamadaAdaptive} likewise use large Cayley
coordinates as a warning signal to change the center and restart the
Euclidean solver on the Stiefel manifold.}
The following estimates motivate a radius-based trigger without evaluating
Jacobian condition numbers.

At a stationary point, the chart differential controls the conditioning
of the pullback Hessian through the following bound.
\begin{proposition}[Hessian conditioning under a chart]\label{prop:conditioning}
Let $x=G(y)$ be stationary, with positive-definite Riemannian Hessian.
Then
\[
 \kappa_2(\nabla^2h(y))
 \le\kappa_2(J_y)^2\,
       \kappa_2(\widehat{\Hess f}(x)).
\]
\end{proposition}
\begin{proof}
At stationarity, Lemma~\ref{lem:deriv} gives
$\nabla^2h=J_y^\top\widehat{\Hess f}J_y$.
The Rayleigh quotient bounds on its extreme eigenvalues give the claim.
\end{proof}

Proposition~\ref{prop:conditioning} motivates controlling Jacobian
conditioning through the chart radius. Away from stationarity, the
additional covariant-derivative term in Lemma~\ref{lem:deriv} must also be
controlled. We next derive radius bounds under the stated geometric
assumptions. Detailed proofs and
coordinate conversions are collected in Appendix~\ref{app:radius-proofs}.
All bounds concern the specified metric and parameter norm.

\Needspace{8\baselineskip}
The radius can be written in terms of the chosen manifold and \rev{coordinate representation}, and sharper bounds are available for particular charts.
For homogeneous Cayley charts, the following transversality assumption
gives a general sufficient bound; the subsequent product and graph-chart
results improve its quartic dependence to a quadratic one.
\rev{The next proposition is stated for homogeneous spaces satisfying
the transversality assumption below. In particular, it applies to the
real Stiefel Cayley chart with the normal homogeneous metric.}
\begingroup
\begin{proposition}[Chart conditioning on Stiefel-type manifolds]\label{thm:homogeneous-quartic}
Let $M=\SO(N)/H$ have the normal homogeneous
Frobenius metric. Set $\mathfrak m=\mathfrak h^\perp$, use orthonormal
coordinates on $\mathfrak m$, and define
$\Phi_x(Y)=g\cdot\pi((I+Y)(I-Y)^{-1})$ for $x=g\cdot\pi(I)$.

With $R_Y=(I-Y)^{-1}$ and horizontal
projection $P_{\mathfrak m}$, assume for all $Y,Z\in\mathfrak m$ that
\begin{equation}\label{eq:cayley-transversality}
 \|R_Y^\top ZR_Y\|_F \leq(1+\|Y\|_2^2)\|P_{\mathfrak m}(R_Y^\top ZR_Y)\|_F.
\end{equation}

Then, uniformly in $x$, for $r=\|Y\|_F$,
\begin{equation}\label{eq:homogeneous-quartic}
 \kappa_2(D\Phi_x(Y))\leq(1+\|Y\|_2^2)^2\leq(1+r^2)^2.
\end{equation}
For $\rho>1$, the sufficient radius is
\begin{equation}\label{eq:homogeneous-radius}
 r\leq r_\rho=\sqrt{\sqrt\rho-1} \quad\Longrightarrow\quad\kappa_2(D\Phi_x(Y))\leq\rho.
\end{equation}
\end{proposition}
\endgroup

The real Stiefel Cayley chart satisfies
\eqref{eq:cayley-transversality} in the normal homogeneous metric.
Table~\ref{tab:additional-radius-bounds} also collects the product,
complex Stiefel, and generalized Stiefel bounds; their proofs and coordinate
conversions are in Appendix~\ref{app:radius-proofs}.
For rotations, $Y_i=[w_i]_\times$ and $r^2=\sum_i\|w_i\|_2^2$;
for oblique manifolds, use columnwise stereographic coordinates with
$r^2=\sum_i\|b_i\|_2^2$. Complex Stiefel uses the normal metric on
$U(N)/U(N-p)$. Whitening $X\mapsto B^{1/2}X$ pulls this convention
back to the generalized Stiefel manifold; an ambient target metric
introduces $K_B=\sqrt{2\kappa_2(B)}$, with $\rho>K_B$ in the radius bound.
The normalized Stiefel radius is $r=\|Y\|_F$; independent-entry coordinates
require the conversions in the appendix. All bounds concern the stated
chart and metric, and away from stationarity the additional term in
Lemma~\ref{lem:deriv} remains relevant.

\Needspace{6\baselineskip}
For Grassmann, the graph chart \eqref{eq:grassmann-chart}
admits a quadratic conditioning bound.
\begingroup
\begin{proposition}[Grassmann manifold]\label{thm:grassmann-graph}
Consider $G_{Q_0}$ in \eqref{eq:grassmann-chart} with the standard metric
of Section~\ref{ssec:test-charts}. For $y=\operatorname{vec}(B)$, write
$J_B=J_y$ and $r=\|B\|_F$. Uniformly over the choice of $Q_0$,
\begin{equation}\label{eq:grassmann-graph-bounds}
\begin{aligned}
 \frac{1}{1+\|B\|_2^2}&\leq\sigma_{\min}(J_B)\leq\sigma_{\max}(J_B)\leq1,\\
 \kappa_2(J_B)&\leq1+\|B\|_2^2\leq1+r^2.
\end{aligned}
\end{equation}
Consequently, for $\rho>1$,
\begin{equation}\label{eq:grassmann-graph-radius}
 r\leq r_\rho=\sqrt{\rho-1} \quad\Longrightarrow\quad\kappa_2(J_B)\leq\rho.
\end{equation}
\end{proposition}
\endgroup

\begin{table}[H]\centering\small
\caption{\rev{Jacobian condition-number bounds and sufficient radii for the stated charts.}}
\label{tab:additional-radius-bounds}
\begin{tabular}{p{0.47\linewidth}ll}\toprule
Geometry and metric & Bound on $\kappa_2(J)$ & Sufficient radius\\\midrule
Real Stiefel manifold, normal & $(1+r^2)^2$ & $\sqrt{\sqrt\rho-1}$\\
$\SO(3)^q$, product Frobenius & $1+r^2$ & $\sqrt{\rho-1}$\\
$\OB(d,q)$, product round & $1+r^2$ & $\sqrt{\rho-1}$\\
Complex Stiefel \rev{manifold}, normal & $(1+r^2)^2$ & $\sqrt{\sqrt\rho-1}$\\
Generalized Stiefel \rev{manifold}, whitened normal & $(1+r^2)^2$ & $\sqrt{\sqrt\rho-1}$\\
Generalized Stiefel \rev{manifold}, ambient & $K_B(1+r^2)^2$ & $\sqrt{\sqrt{\rho/K_B}-1}$\\
Grassmann, graph; standard & $1+r^2$ & $\sqrt{\rho-1}$\\
\bottomrule\end{tabular}
\end{table}

The restart test uses the coordinate displacement
\begin{equation}\label{eq:bfgs-chart-radius}
 r_k(y)=\|y-y_0\|_2,\qquad y_0=0,\qquad R=3.
\end{equation}
The common threshold was selected during preliminary training. Checks
after every ten accepted steps may allow overshoot, so this is a restart
safeguard rather than a hard constraint on all iterates.

%
After an inner exit, test the Riemannian gradient at
$x_{k+1}=G_{x_k}(y)$. If it is not small enough, rebuild the chart and
clear the history. Retain only the positive scalar scale, with the
clipped adjustment
\[
 \gamma_{\rm new}\approx\gamma_{\rm old}
 \left(\frac{\|g_{\rm old}\|_2}{\|g_{\rm new}\|_2}\right)^2.
\]
The two gradients are evaluated at the same manifold point in the old
and new charts. This approximately preserves the predicted initial
decrease $\gamma\|g\|_2^2$; anisotropic chart changes make it a heuristic
scale correction, not a curvature transformation. The scale is clipped
to remain positive and finite, with default value one when unavailable.

Algorithm~\ref{alg:restart-summary} summarizes the chart-restart procedure.
The inner stopping test ends a chart; convergence is certified by the
Riemannian gradient stopping criterion.
\begin{algorithm}[tbp]
\caption{Chart Restart BFGS}\label{alg:restart-summary}
\begin{algorithmic}[1]
\Require $x_0\in\M$, tolerance $\tau>0$, budgets $B,J$, fixed radius $R>0$
\State $x\gets x_0$; initialize a positive scalar scale $\gamma$
\For{$k=0,\ldots,B-1$}
\If{$\norm{\grad f(x)}<\tau$}
\State \Return $x$ (converged)
\EndIf
\State Construct $G_x$; set $h=f\circ G_x$, $y=0$; clear curvature history
\State Calibrate the coordinate tolerance by \eqref{eq:bfgs-dynamic-gtol}
\State Run Euclidean BFGS (L-BFGS in the experiments), initialized with scale $\gamma$
\State\label{line:restart-radius} \rev{At scheduled checks, stop} the inner solve \rev{if $\norm{y}_2\ge R$ (radius exit)}
\State \rev{Also stop} on \rev{the inner} stopping test, failure, or $J$ accepted steps
\State $x^+\gets G_x(y)$ using the last valid inner iterate
\State If $x^+$ passes the Riemannian gradient stopping criterion above, \Return $x^+$ (converged)
\State If recovery fails with $y=0$, \Return $x$ (stalled)
\State $x\gets x^+$; retain only the clipped positive scalar scale $\gamma$ for the next chart
\EndFor
\State \Return $x$ (budget exhausted)
\end{algorithmic}
\end{algorithm}

\subsection{Convergence Analysis}
\rev{We establish stationarity of accumulation points of the
chart centers. We then show that radius restarts stop near a local
minimizer and derive local superlinear convergence of full-memory BFGS
in the final chart.}

Let $x_k$ be the center of the $k$th chart, and write
\begin{equation}\label{eq:restart-notation}
 h_k=f\circ G_{x_k},\qquad y_{k,0}=0,\qquad
 x_{k,j}=G_{x_k}(y_{k,j}).
\end{equation}
Here $k$ counts chart constructions and $j$ counts accepted steps within
a chart. If a radius restart occurs after $j_k$ accepted steps, then
$x_{k+1}=x_{k,j_k}$. If only finitely many restarts occur, the last chart
may contain infinitely many inner iterations.

We \rev{analyze} the infinite-horizon version of
Algorithm~\ref{alg:restart-summary}, with $J=\infty$, no total chart
budget, and only radius-triggered chart changes. Coordinate tolerances
are tightened within the same chart when necessary. Each chart
satisfies $G_x(0)=x$ and starts with \rev{$H_{k,0}=\gamma_{k,0}I$}, where
$0<\underline\gamma\le\gamma_{k,0}\le\bar\gamma$.
Line searches succeed and satisfy the \rev{Armijo--}Wolfe conditions with
fixed $0<c_1<1/2$ and $c_1<c_2<1$. \rev{Accepted steps obey}
$\alpha_{k,j}\le\bar\alpha$\rev{, where $\bar\alpha\ge1$; for the
superlinear rate, the search tries} $\alpha=1$ \rev{first} and \rev{accepts it
whenever it satisfies} these \rev{conditions}. These \rev{successful bounded
searches are theoretical requirements, separate} from the finite budgets
and failure exits in the \rev{implementation}.

\begin{assumption}\label{ass:restart-global}
\rev{The objective satisfies $f\in C^1(\M)$.}
\begin{enumerate}[label=(\arabic*)]
\item The initial sublevel set
$S=\{x\in\M:f(x)\le f(x_0)\}$ is compact.
\item The chart family $(x,y)\mapsto G_x(y)$ is jointly $C^2$ \rev{near}
$S\times\overline B_T(0)$, \rev{with smooth local frame choices. The} parameter
ball contains every first accepted line-search segment.
\end{enumerate}
\end{assumption}
\begin{remark}
\rev{Part~(1) is automatic when $\M$ is compact. On a noncompact
manifold, it suffices that $f$ have compact sublevel sets, for example
when it is coercive. This is the usual level-set boundedness requirement
in descent analysis}~\cite[Chapter 3]{NW2006}\rev{.}

\rev{For part~(2), the Cayley and graph formulas, as well as the
stereographic formula, are jointly smooth near
$S\times\overline B_T(0)$ for every finite $T$, with smooth local
choices of center frames. Thus their regularity imposes no restriction
on the choice of a finite $T$. To ensure that the ball also contains
every first accepted line-search segment, choose}
\revdisplay{\[
 T\ge\bar\alpha\bar\gamma G_0,
 \qquad G_0:=\max_{x\in S}\|\nabla(f\circ G_x)(0)\|_2<\infty.
\]}
\rev{Compactness of $S$ and continuity at the chart origins give
$G_0<\infty$. Since the first direction is
$-\gamma_{k,0}\nabla h_k(0)$, the scale and step-length bounds imply
$\|y_{k,1}\|_2\le\bar\alpha\bar\gamma G_0\le T$.
The entire segment from $0$ to $y_{k,1}$ therefore lies in
$\overline B_T(0)$, verifying part~(2).}
\end{remark}

\rev{The next theorem shows that, if radius restarts occur infinitely
often, the gradients at the chart centers vanish and
stationary accumulation points exist.}
\begin{theorem}[Stationary accumulation points of restart centers]\label{thm:restart-stationary}
Under Assumption~\ref{ass:restart-global}, an infinite restart sequence
satisfies $\norm{\grad f(x_k)}\to0$ and has \rev{at least one}
stationary accumulation point.
\end{theorem}
\begin{proof}
\rev{Armijo decrease and recentering preserve membership in} $S$\rev{, so} $f$ \rev{is
bounded below along the iteration. Joint regularity, compactness, and
invertibility of the chart differentials give, using a finite cover
by smooth local frame choices,}
\[
 m_0=\min_{x\in S}
 \sigma_{\min}(DG_x(0))>0\rev{.}
\]
The \rev{pullback gradient identity yields
$\norm{g_{k,0}}_2\ge m_0\norm{\grad f(x_k)}$, where}
$g_{k,0}=\nabla h_k(0)$\rev{.}

\rev{Fix} $\varepsilon>0$\rev{. Since $f\in C^1(\M)$, the pullback gradients
are uniformly continuous on the compact center and parameter sets
in Assumption}~\ref{ass:restart-global}\rev{. Thus there is a
$\delta_\varepsilon>0$ such that, uniformly over admissible centers,}
\revdisplay{\[
 \norm{\nabla h_x(y)-\nabla h_x(0)}_2
 <(1-c_2)m_0\varepsilon
 \quad\text{whenever }\norm y_2<\delta_\varepsilon
 \text{ and }\norm y_2\le T.
\]}
\rev{For the first direction $p_{k,0}=-\gamma_{k,0}g_{k,0}$, the} Wolfe
curvature condition \rev{implies}
\revdisplay{\[
 \norm{\nabla h_k(y_{k,1})-g_{k,0}}_2
 \ge(1-c_2)\norm{g_{k,0}}_2.
\]}
\rev{Consequently, if $\norm{\grad f(x_k)}\ge\varepsilon$, then
$\norm{y_{k,1}}_2\ge\delta_\varepsilon$}.
Armijo decrease and subsequent monotonicity \rev{give}
\revdisplay{\begin{equation}\label{eq:restart-outerdecrease}
 f(x_k)-f(x_{k+1})
 \ge c_1\norm{y_{k,1}}_2\norm{g_{k,0}}_2
 \ge c_1\delta_\varepsilon m_0\varepsilon>0.
\end{equation}}
This \rev{uniform decrease can occur only finitely often because} $f$ \rev{is
bounded below}. \rev{Since $\varepsilon$ was arbitrary,}
$\norm{\grad f(x_k)}\to0$. Compactness supplies a convergent
subsequence, and gradient continuity \rev{makes its limit stationary.}\hfill\proofbox\NoEndMark

\end{proof}

\rev{We next} impose local conditions \rev{that prevent a chart centered near a
minimizer from leaving its regular neighborhood}.
\begin{assumption}\label{ass:restart-local}
\rev{The objective} $f$ \rev{is twice continuously differentiable on a neighborhood
of} $x^\star$\rev{.} The point $x^\star$ is stationary and satisfies
\begin{equation}\label{eq:restart-nondegenerate}
 \Hess f(x^\star)\succ0.
\end{equation}
The pullback Hessians are locally Lipschitz near $(x^\star,0)$,
with uniform constants for nearby centers.
\end{assumption}
\begin{remark}
Assumption~\ref{ass:restart-local} specifies a nondegenerate local minimizer
and locally Lipschitz pullback Hessians, as in standard local BFGS and
quasi-Newton analyses~\cite{BDM1973}\cite[Section 6.4]{NW2006}.
The Lipschitz condition holds, for example, when the objective and centered
chart maps are locally $C^3$.
\end{remark}

\begin{lemma}[Local \rev{stability} and \rev{no further} radius restarts]\label{thm:restart-nonescape}
Suppose Assumptions~\ref{ass:restart-global} and~\ref{ass:restart-local}
hold. Each new chart uses full-memory BFGS initialized at $0$ with empty
history and $H_{x,0}=\gamma_xI$, under the scale and successful bounded
Armijo--Wolfe conventions above. For \rev{every} fixed $R>0$, there is a
neighborhood $U_0$ of $x^\star$ such that a chart constructed at
$x\in U_0$ satisfies
\begin{equation}\label{eq:restart-stability}
 \sup_{j\ge0}\norm{y_{x,j}-y_x^\star}_2
 \le\sqrt{L/\mu}\,\norm{y_x^\star}_2,
 \qquad y_{x,j}\to y_x^\star\rev{,}
\end{equation}
\rev{where $G_x(y_x^\star)=x^\star$ and} $0<\mu\le L$ \rev{are uniform local
Hessian bounds.} In particular, $\norm{y_{x,j}}_2<R$ for every \rev{accepted}
inner iterate, so no radius restart occurs.
\end{lemma}
\begin{proof}
\rev{Define $\Phi(x,y)=G_{x^\star}^{-1}(G_x(y))$ near $(x^\star,0)$,
using the local inverse of $G_{x^\star}$ near $x^\star$. Since
$\Phi(x^\star,0)=0$ and
$D_y\Phi(x^\star,0)=(DG_{x^\star}(0))^{-1}DG_{x^\star}(0)=I$, applying
the implicit function theorem to $\Phi$~\cite[Theorem C.40]{Lee2013}
gives local coordinates} $y_x^\star$ \rev{with $G_x(y_x^\star)=x^\star$ and
$y_x^\star\to0$ as $x\to x^\star$. By Lemma}~\ref{lem:deriv}
\rev{and Assumption}~\ref{ass:restart-local}\rev{, there exist uniform
constants} $r>0$ \rev{and} $0<\mu\le L$ \rev{such that}
\begin{equation}\label{eq:restart-localconvex}
 \mu I\preceq\nabla^2h_x(y)\preceq LI,
 \qquad \norm{y-y_x^\star}_2\le2r,
\end{equation}
\rev{for all centers $x$ sufficiently close to $x^\star$.}
Lemma~\ref{lem:restart-extension} \rev{in Appendix}~\ref{app:restart-stability}
\rev{constructs strongly convex} extensions $\widetilde h_x=h_x$ on
$B_r(y_x^\star)$ and \rev{proves} a uniform bound
$\widetilde H_{x,j}\preceq M_HI$ for their full-memory BFGS \rev{iterates.
This extension is only an analytical device; the algorithm evaluates} $h_x$.

Put $C=\sqrt{L/\mu}$ and $\delta_x=\norm{y_x^\star}_2$.
\rev{Choose} $U_0$ so that
\[
 (1+\bar\alpha M_HL)C\delta_x<r,
 \qquad (1+C)\delta_x<R.
\]
Suppose an accepted step is the first to leave $B_r(y_x^\star)$.
\rev{Every} preceding \rev{step}, including \rev{its} endpoint gradients, \rev{is a} valid
Armijo--Wolfe \rev{step} for $\widetilde h_x$ \rev{and} can be continued as an
auxiliary BFGS sequence\rev{. Hence} $H_{x,j}\preceq M_HI$ \rev{at the} current
\rev{iterate}. Descent and \eqref{eq:restart-localconvex} give
$\norm{y_{x,j}-y_x^\star}_2\le C\delta_x$\rev{, and therefore}
\[
 \norm{y_{x,j+1}-y_x^\star}_2
 \le(1+\bar\alpha M_HL)C\delta_x<r,
\]
\rev{a contradiction}. \rev{All} accepted \rev{points stay} in the ball, \rev{where} descent
gives \eqref{eq:restart-stability}. \rev{The} convergence \rev{estimate} in
Lemma~\ref{lem:restart-extension} also gives $y_{x,j}\to y_x^\star$
\rev{and}
\revdisplay{\[
 \sum_{j=0}^\infty\norm{y_{x,j}-y_x^\star}_2<\infty.
\]}
Finally $\norm{y_{x,j}}_2\le(1+C)\delta_x<R$.\hfill\proofbox\NoEndMark

\end{proof}

The \rev{next lemma records} the \rev{Euclidean rate used after} local \rev{stability
has been established}. Its \rev{summability hypothesis is supplied by the
preceding proof, rather than by} a \rev{separate compact-iterate assumption}.
\begin{lemma}[Local $q$-superlinear convergence in a fixed chart]\label{thm:super}
\rev{Let} $h=f\circ G$ \rev{in a fixed chart.} Suppose full-memory BFGS iterates
converge to $y^\star$, with
$\nabla^2h(y^\star)\succ0$\rev{, locally} Lipschitz \rev{Hessian}, and
\rev{$\sum_j\norm{y_j-y^\star}_2<\infty$. Use an Armijo--Wolfe line search
with} $0<c_1<1/2$\rev{,} $c_1<c_2<1$\rev{, which tries and accepts the unit step
whenever it is admissible. Then either the method terminates at}
$y^\star$\rev{, or} the unit step is eventually accepted \rev{and}
\revdisplay{\[
 \norm{y_{j+1}-y^\star}_2=o(\norm{y_j-y^\star}_2).
\]}
\rev{The same conclusion holds for a cautious variant once all updates are
accepted}.
\end{lemma}
\begin{proof}
The \rev{pairs} $s_j=y_{j+1}-y_j$ \rev{and
$z_j=\nabla h(y_{j+1})-\nabla h(y_j)$ are exact Euclidean} secant
pairs\rev{. The Hessian regularity} and \rev{summability give} the \rev{BFGS}
Dennis--Mor\'e \rev{condition}~\cite[Theorem 6.6]{NW2006}\rev{,}
\revdisplay{\[
 \frac{\norm{(H_j^{-1}-\nabla^2h(y^\star))p_j}_2}
      {\norm{p_j}_2}\longrightarrow0,
 \qquad p_j=-H_j\nabla h(y_j).
\]}
\rev{The unit step then satisfies} the \rev{Wolfe conditions for all sufficiently
large} $j$~\cite[Theorem 3.6]{NW2006}\rev{;} the \rev{stated search convention makes
it the accepted step. The} Dennis--Mor\'e characterization \rev{yields the
claimed rate}~\cite{DennisMore1974}\rev{.}
\end{proof}

\rev{Together, the two lemmas show that the full-memory version of
Algorithm~\ref{alg:restart-summary} needs no further radius restarts
and converges $q$-superlinearly once a chart is centered sufficiently
near $x^\star$.}
\begin{theorem}[\rev{No further radius restarts and local superlinear convergence}]\label{prop:restart-finite}
\rev{Consider the infinite-horizon full-memory version of
Algorithm}~\ref{alg:restart-summary} \rev{under
Assumptions}~\ref{ass:restart-global} \rev{and}~\ref{ass:restart-local}\rev{.
There} is \rev{a neighborhood} $U_0$ \rev{of} $x^\star$ \rev{such that no radius restart
occurs in any chart constructed at a center $x_k\in U_0$.
The iterates in that chart either terminate at} $x^\star$ \rev{or converge}
$q$-\rev{superlinearly to} $x^\star$ \rev{in Riemannian distance.}
\end{theorem}
\begin{proof}
\rev{For a center in} $U_0$\rev{, Lemma}~\ref{thm:restart-nonescape} \rev{gives
$\norm{y_{k,j}}_2<R$ for every accepted iterate, so the radius test
never triggers. That lemma also establishes} convergence \rev{and summability
of} the \rev{coordinate errors}, with \rev{all cautious updates accepted locally}.
\rev{Lemma}~\ref{thm:super} \rev{therefore applies}. A regular chart is locally
bi-Lipschitz between Euclidean and Riemannian distances, so
\revdisplay{\[
 \frac{\dist(G_{x_k}(y_{k,j+1}),x^\star)}
      {\dist(G_{x_k}(y_{k,j}),x^\star)}\longrightarrow0.
\]}
\end{proof}

\begin{remark}\label{prop:stall}
In a fixed chart, the Euclidean secant pairs require no vector transport.
Consequently, the local full-memory BFGS analysis follows directly from
Euclidean theory under the hypotheses of Lemma~\ref{thm:super}, without
imposing isometry or locking conditions on a vector transport.
\end{remark}

\section{Numerical experiments}\label{sec:num}

The test suite comprises 15 problem families E01--E15 adapted from the
Manopt and Pymanopt example collections~\cite{Manopt2014,Pymanopt2016}.
They include subspace and eigenvalue problems, matrix factorizations,
sphere packing, semidefinite optimization, statistical objectives, and
low-rank recovery, and so on. The complete list appears in
Table~\ref{tab:bfgs-multiseed-problems}. We compare chart-restarted Euclidean L-BFGS
(Euc-BFGS; Algorithm~\ref{alg:restart-summary}) with Riemannian
limited-memory BFGS (RiemBFGS) and Riemannian steepest descent (RiemSD).
The two baselines are Python implementations adapted from Manopt's
RLBFGS and steepest-descent solvers~\cite{Manopt2014}; the coordinate
\rev{L-BFGS recursion follows the formulation in~\cite{LiuNocedal1989}.}
All methods are evaluated in the same NumPy/SciPy environment with
gradients. The ambient objective gradient and the coordinate
pullback gradient are both computed from closed-form expressions. For an
ambient extension $\bar f$ of the objective, the latter applies the
explicit adjoint differential,
$\nabla h(y)=DG_x(y)^*\nabla\bar f(G_x(y))$.
The ambient and pullback gradients pass directional-derivative checks on all
150 instances. Computations use one CPU thread for the underlying
numerical libraries.

\subsection{Experimental setup}\label{sec:bfgs-multiseed}
For each problem family, ten random seeds are used to generate the
matrices or vectors defining the problem. Different random seeds yield
different problem instances, giving a total of 150 cases. All three
methods use the same problem data and starting point for each case.
A solve is successful when
\[
 \|\grad f(x)\|<\tau,\qquad \tau=10^{-6}.
\]
Each solver is
warmed up before three timed repeats on each instance; the median is the
reported solve time. Method order is randomized, yielding 1,350 measured
solves. Euclidean L-BFGS uses $B=20$ chart segments with at most $J=100$
accepted steps each, including manifolds \rev{with coordinates defined on all of} $\R^d$; RiemBFGS and
RiemSD each allow at most 2,000 iterations. The coordinate radius is
checked every ten accepted steps, without an early-step exception.
The radius, selected during preliminary training and fixed throughout
the reported experiments, is $3$ for all top-level manifolds.

Unlike Algorithm~\ref{alg:fixed-bfgs}, the experimental solver does not
store a dense inverse-Hessian matrix. It retains only $\ell$ coordinate
secant pairs and applies their BFGS updates through the two-loop recursion.
Algorithm~\ref{alg:lbfgs-inner} specifies the inner iteration called in
Algorithm~\ref{alg:restart-summary}. The stopping test, cautious update,
and radius check are shown explicitly.
\begin{algorithm}[tbp]
\caption{Limited-memory BFGS within one chart}\label{alg:lbfgs-inner}
\begin{algorithmic}[1]
\Require $h=f\circ G_x$, tolerance $\eta$, memory $\ell$, scale $\gamma>0$, budget $J$, radius $R$, cautious constant $c>0$
\State $y\gets0$; $g\gets\nabla h(0)$; $\mathcal H\gets\varnothing$
\For{$j=0,\ldots,J-1$}
\State\label{line:inner-stop} If $\|g\|_2\leq\eta$, \Return $y$ (inner tolerance)
\State\label{line:inner-direction} $p\gets-\operatorname{TwoLoop}(g;\mathcal H,\gamma)$; if descent fails, clear history and use $p=-g$
\State Obtain a safeguarded Wolfe step $y^+=y+\alpha p$; on failure, \Return $y$
\State $g^+\gets\nabla h(y^+)$; $s\gets y^+-y$; $z\gets g^+-g$
\If{$s^\top z>0$ and $s^\top z/(s^\top s)\geq c\|g\|_2$}
\State\label{line:inner-curvature} Store $(s,z)$; retain the newest $\ell$ pairs and update the positive scale
\EndIf
\State $y\gets y^+$; $g\gets g^+$
\State\label{line:inner-radius} At scheduled checks, \Return $y$ if the radius restart test is satisfied
\EndFor
\State \Return $y$ (inner budget)
\end{algorithmic}
\end{algorithm}

\rev{At the origin of each chart, we calibrate the coordinate tolerance by}
\revdisplay{\begin{equation}\label{eq:bfgs-dynamic-gtol}
 q_k=\frac{\norm{\nabla h_k(0)}_2}{\norm{\grad f(x_k)}},
 \qquad \eta_k\approx\tau q_k.
\end{equation}}
\rev{If this accepts the origin before the Riemannian gradient stopping
criterion holds, we tighten $\eta_k$ and retry.}
Algorithm~\ref{alg:lbfgs-inner},
line~\ref{line:inner-stop}, tests the Euclidean norm of the coordinate
gradient using the calibrated tolerance of
\eqref{eq:bfgs-dynamic-gtol}; the outer test above certifies manifold
convergence. Line~\ref{line:inner-direction} replaces the dense BFGS
matrix update in Algorithm~\ref{alg:fixed-bfgs} by a two-loop application.
Line~\ref{line:inner-curvature} implements the Li--Fukushima cautious
update~\cite{LiFukushima2001}. If the cautious test
fails, the new pair is not stored. The separate positive-curvature check
is retained as a numerical safeguard.
Line~\ref{line:inner-radius} checks the coordinate radius periodically.
The memory length is $\ell=30$, and the run budgets are those specified
above. A non-descent direction is replaced by the negative coordinate
gradient, and inner-solver failure returns the last valid point to the
outer stationarity test. For SPD variables, an orthonormal symmetric
coordinate basis ensures $\|S(y)\|_F=\|y\|_2$. All cases in the stated suite, including unsuccessful solves, are
included in the comparisons below.

\subsection{Results and performance profiles}\leavevmode\par
\begin{table}[H]\centering\small
\caption{\rev{Success counts and geometric-mean solve times. Bold indicates the
largest success count or smallest time in each column.
Time differences are computed before rounding.}}
\label{tab:bfgs-multiseed}
\setlength{\tabcolsep}{4pt}
\begin{tabular}{lrrr}
\toprule
Method & \shortstack{Successful cases\\(total 150)} & \shortstack{Geometric-mean time (ms)\\Own successes} & \shortstack{Geometric-mean time (ms)\\All-method successes (115)}\\
\midrule
Euc-BFGS & \rev{\textbf{150}} & \rev{\textbf{28.4}} & \rev{\textbf{23.3}}\\
RiemBFGS & 132 & 58.4 & 52.1\\
RiemSD & 117 & 49.6 & 49.3\\
\bottomrule
\end{tabular}

\end{table}

\begin{table}[tbp]\centering\small
\caption{\rev{Results for E01--E15. Entries give successful cases out of ten;
geometric-mean time (ms) over successful cases. Bold marks the best
values in each row; a dash indicates no success.}}
\label{tab:bfgs-multiseed-problems}
\setlength{\tabcolsep}{4pt}
\begin{tabular}{rlrrr}
\toprule
ID & Problem & Euc-BFGS & RiemBFGS & RiemSD\\
\midrule
E01 & Dominant subspace & \rev{\textbf{10}}; \rev{\textbf{41.1}} & \rev{\textbf{10}}; 41.2 & \rev{\textbf{10}}; 54.2\\
E02 & Complex dominant subspace & \rev{\textbf{10}}; 55.1 & \rev{\textbf{10}}; \rev{\textbf{45.2}} & 6; 49.5\\
E03 & Generalized eigenproblem & \rev{\textbf{10}}; \rev{\textbf{32.9}} & 0; --- & 0; ---\\
E04 & Nonlinear eigenspace & \rev{\textbf{10}}; 50.3 & \rev{\textbf{10}}; 42.6 & \rev{\textbf{10}}; \rev{\textbf{35.7}}\\
E05 & Truncated SVD & \rev{\textbf{10}}; 56.4 & 3; 144.2 & 1; \rev{\textbf{39.9}}\\
E06 & Sparse PCA & \rev{\textbf{10}}; \rev{\textbf{10.0}} & \rev{\textbf{10}}; 11.2 & \rev{\textbf{10}}; 10.3\\
E07 & Thomson problem & \rev{\textbf{10}}; \rev{\textbf{26.3}} & \rev{\textbf{10}}; 53.4 & \rev{\textbf{10}}; 43.8\\
E08 & Sphere packing & \rev{\textbf{10}}; \rev{\textbf{117.9}} & \rev{\textbf{10}}; 224.0 & 5; 450.6\\
E09 & Elliptope SDP & \rev{\textbf{10}}; \rev{\textbf{56.3}} & \rev{\textbf{10}}; 111.1 & 6; 113.7\\
E10 & Max-cut & \rev{\textbf{10}}; \rev{\textbf{21.2}} & \rev{\textbf{10}}; 49.2 & \rev{\textbf{10}}; 33.6\\
E11 & Generalized Procrustes & \rev{\textbf{10}}; \rev{\textbf{16.5}} & \rev{\textbf{10}}; 659.9 & \rev{\textbf{10}}; 328.1\\
E12 & Smoothed robust PCA & \rev{\textbf{10}}; \rev{\textbf{16.2}} & \rev{\textbf{10}}; 21.1 & \rev{\textbf{10}}; 18.4\\
E13 & SPD intrinsic mean & \rev{\textbf{10}}; 6.9 & \rev{\textbf{10}}; \rev{\textbf{6.2}} & \rev{\textbf{10}}; 21.4\\
E14 & Essential SVD & \rev{\textbf{10}}; \rev{\textbf{7.3}} & \rev{\textbf{10}}; 73.0 & \rev{\textbf{10}}; 52.4\\
E15 & Low-rank completion & \rev{\textbf{10}}; \rev{\textbf{60.0}} & 9; 218.6 & 9; 98.9\\
\midrule
Total & & \rev{\textbf{150}}; \rev{\textbf{28.4}} & 132; 58.4 & 117; 49.6\\
\bottomrule
\end{tabular}

\end{table}

Tables~\ref{tab:bfgs-multiseed} and~\ref{tab:bfgs-multiseed-problems}
show that Euc-BFGS solves all 10 instances in each family of the stated
suite. RiemBFGS and RiemSD solve 132 and 117 instances, respectively.
Euc-BFGS has a clear robustness advantage on E03 (generalized eigenproblems)
and E05 (truncated SVD).
For a comparison on identical instances, we restrict the timing analysis
to the 115 cases solved by all three methods. Their geometric-mean solve
times are $23.3$~ms for Euc-BFGS, $52.1$~ms for RiemBFGS, and $49.3$~ms
for RiemSD. Thus Euc-BFGS \rev{uses approximately $45\%$ and $47\%$ of}
the geometric-mean solve \rev{times of} RiemBFGS and RiemSD\rev{, respectively,}
on this common set\rev{---roughly half the time of either baseline}. Means over each method's own successes are also reported
in Table~\ref{tab:bfgs-multiseed} for completeness.

\begingroup
\setlength{\abovedisplayskip}{6pt}
\setlength{\belowdisplayskip}{6pt}
\setlength{\intextsep}{6pt}
\setlength{\textfloatsep}{12pt}
\setlength{\abovecaptionskip}{4pt}
\Needspace{6\baselineskip}
\rev{To compare the distribution of outcomes, we use performance profiles
as defined in~\cite{DolanMore2002}.} For case $i$ and solver
$s$, let $m_{i,s}$ be the performance measure and set
\[
 r_{i,s}=\frac{m_{i,s}}{\min_t m_{i,t}},\qquad
 P_s(a)=\frac{1}{150}\#\{i:r_{i,s}\leq2^a\},\qquad a\geq0.
\]
For the accuracy comparison, use the terminal normalized residual
\[
 e_{i,s}=\frac{\|\grad f(x_{i,s})\|}{\tau}.
\]
This measures stationarity, rather than error relative to an unknown
optimal objective value. For the time comparison,
$m_{i,s}$ is the median solve time if the solver succeeds and is infinite
otherwise. Cases on which every solver fails contribute zero to every
time profile. A higher curve indicates better relative performance;
the time-profile plateau records successful coverage of the full suite.

\begin{figure}[htbp]\centering
\includegraphics[width=\linewidth]{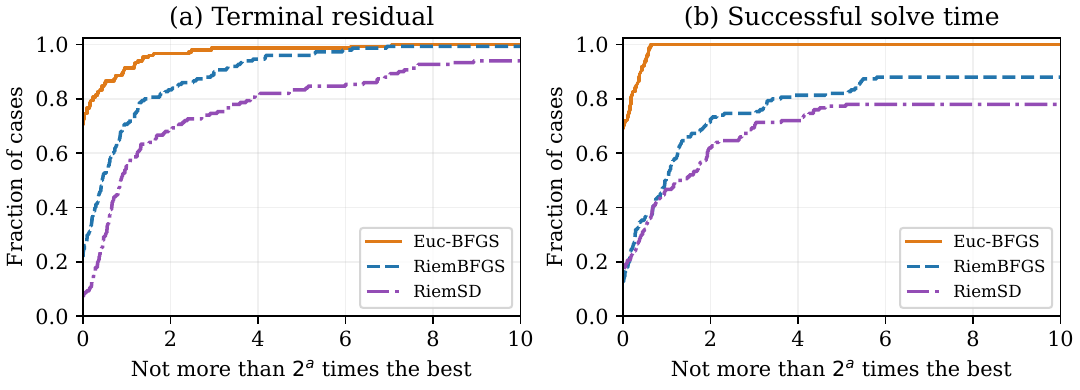}
\caption{Performance profiles over 150 cases: (a) terminal stationarity residual;
(b) successful solve time.}
\label{fig:performance-profiles}
\end{figure}

Figure~\ref{fig:performance-profiles} \rev{shows that Euc-BFGS leads both}
the \rev{terminal-residual and solve-time profiles. It achieves the smallest}
terminal \rev{residual on $106$ cases} and \rev{the fastest} successful solve \rev{on
$104$ cases, while solving all $150$ cases within twice} the \rev{best
successful time. RiemBFGS has a larger pooled geometric-mean time than
RiemSD but a higher time profile, partly because it solves more cases}.

Figure~\ref{fig:three-way-time-seed6} shows \rev{representative convergence}
trajectories for all 15 \rev{problem} families. \rev{Euc-BFGS reaches} the \rev{prescribed
stationarity tolerance in every displayed case}. The Euclidean inner solver uses
$\nabla h(y)$\rev{; its plotted} curve \rev{connects} the recorded \rev{norms of}
$\grad f(G_x(y))$ by straight segments\rev{, without additional per-step}
Riemannian gradient \rev{evaluations}. These segments do not \rev{determine} an
asymptotic convergence order.

\par\endgroup

\begin{figure}[htbp]\centering
\includegraphics[width=\linewidth,height=0.72\textheight,keepaspectratio,trim=0 0 0 30,clip]{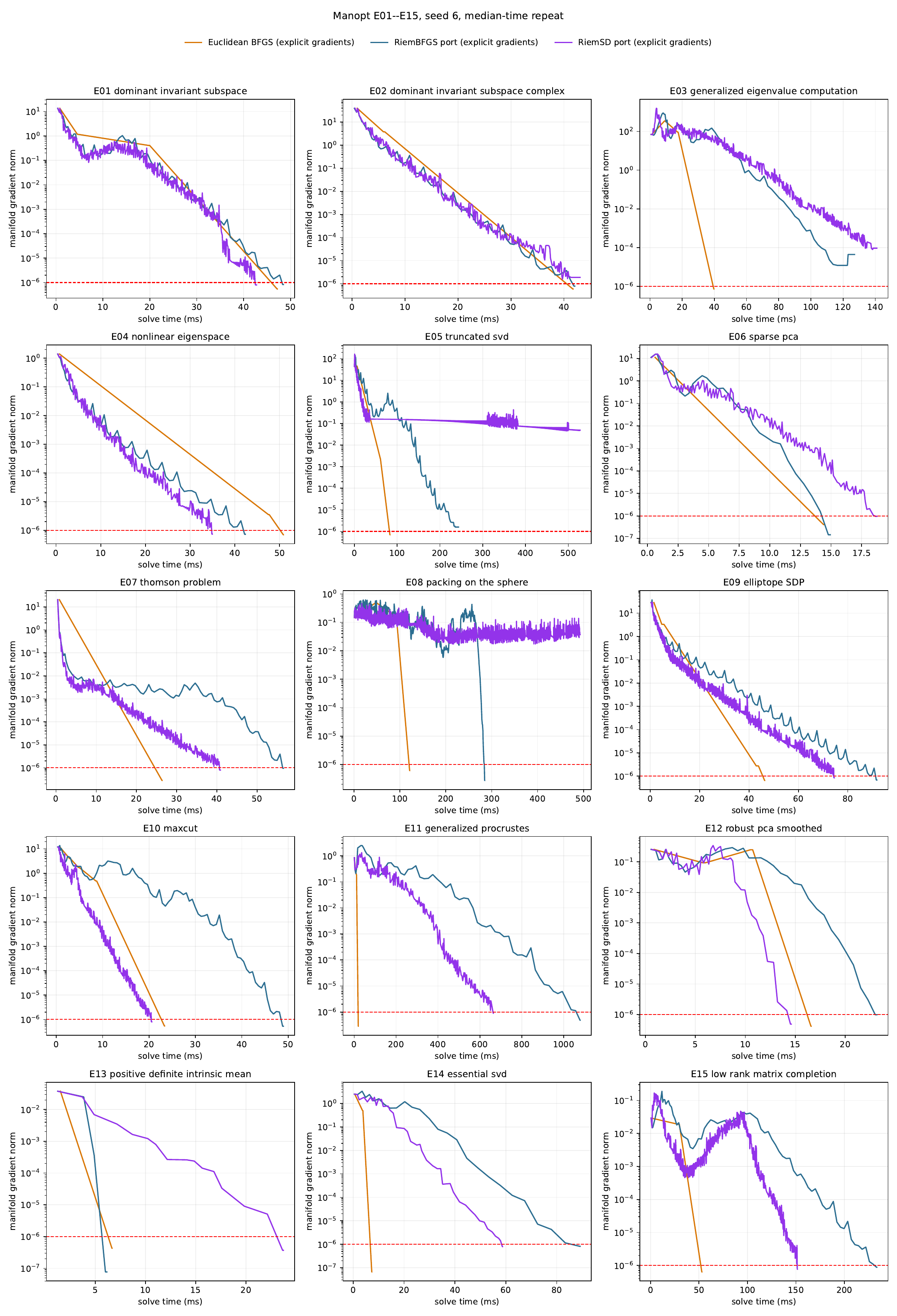}
\setlength{\abovecaptionskip}{5pt}
\caption{Representative convergence trajectories for E01--E15.
The red dashed line marks the stopping tolerance $10^{-6}$.}
\label{fig:three-way-time-seed6}
\end{figure}
\clearpage

\section{Conclusion}\label{sec:conc}
A chart reformulation yields a Euclidean quasi-Newton method using objective values
and first derivatives. Curvature is learned from exact within-chart secant
pairs, and no per-step tangent-vector transport is constructed. The analysis
\rev{establishes stationary-point and minimizer correspondence.
For an infinite restart sequence, we prove that the Riemannian gradients
at the chart centers converge to zero and stationary accumulation
points exist. Local $q$-superlinear convergence is proved for full-memory BFGS
in the final fixed chart. Since the secant pairs are formed in a single
Euclidean space, the analysis requires no vector transport or locking
condition.}
The experiments support the effectiveness of the proposed method across
the tested matrix-manifold problems.

\appendix
\section{Proofs and coordinate details for radius bounds}\label{app:radius-proofs}
The Cayley map satisfies $D\Phi(0)=2D\pi_I|_{\mathfrak m}$ and has zero
covariant acceleration along $t\mapsto\Phi(tZ)$ at the origin. Its
singular values satisfy
\begin{equation}\label{eq:homogeneous-singular-bounds}
 \frac{2}{(1+\|Y\|_2^2)^2}\leq\sigma_{\min}(D\Phi(Y))
 \leq\sigma_{\max}(D\Phi(Y))\leq2.
\end{equation}
\begin{proof}[Proof of \rev{Proposition}~\ref{thm:homogeneous-quartic}]
The expansion $Q(tZ)=I+2tZ+2t^2Z^2+O(t^3)=\exp(2tZ)+O(t^3)$,
together with the normal homogeneous geodesic
$t\mapsto\pi(\exp(2tZ))$ and the isometry $D\pi_I|_{\mathfrak m}$,
gives the assertions at the origin. Fix $Y$, and write $R=R_Y$ and
$k=k_Y$. Since $Q=2R-I$, we have $DQ(Y)[Z]=2RZR$ and
$Q^\top DQ(Y)[Z]=2T$, where $T=R^\top ZR$. The Riemannian submersion
property gives
\begin{equation}\label{eq:quotient-differential}
 \|D\Phi(Y)[Z]\|_g=2\|P_{\mathfrak m}T\|_F.
\end{equation}
Skew-symmetry yields $\|R\|_2\leq1$ and
$\sigma_{\min}(R)=(1+k^2)^{-1/2}$, hence
$\|Z\|_F/(1+k^2)\leq\|T\|_F\leq\|Z\|_F$.
Combining this with \eqref{eq:cayley-transversality} bounds
$\|D\Phi(Y)[Z]\|_g$ between $2\|Z\|_F/(1+k^2)^2$ and $2\|Z\|_F$,
proving \eqref{eq:homogeneous-singular-bounds} and
\eqref{eq:homogeneous-quartic}. Isometric translation preserves the bounds
at other centers, and solving $(1+r^2)^2\leq\rho$ gives
\eqref{eq:homogeneous-radius}.
\end{proof}

\begin{proof}[Verification of the real Stiefel bound in Table~\ref{tab:additional-radius-bounds}]
In standard block coordinates, write
$K=-Y=\left[\begin{smallmatrix}A&B^\top\\-B&0\end{smallmatrix}\right]$,
$R=(I+K)^{-1}$, and
$T=R^\top ZR=\left[\begin{smallmatrix}U&V\\-V^\top&W\end{smallmatrix}\right]$.
For $Z\in\mathfrak m$, the lower-right block of
$Z=(I+K)^\top T(I+K)$ vanishes, giving
$W=-BUB^\top-BV+V^\top B^\top$. With $b=\|B\|_2$, Cauchy--Schwarz gives
\[
 \|W\|_F^2\leq(b^2\|U\|_F+2b\|V\|_F)^2
 \leq(b^4+2b^2)\|P_{\mathfrak m}T\|_F^2.
\]
Thus $\|T\|_F^2=\|P_{\mathfrak m}T\|_F^2+\|W\|_F^2
\leq(1+b^2)^2\|P_{\mathfrak m}T\|_F^2$.
Since $b\leq\|K\|_2=k_Y$, \eqref{eq:cayley-transversality} follows.
\end{proof}

\begin{proof}[Proof of the product bounds in Table~\ref{tab:additional-radius-bounds}]
For a rotation factor, put $R_i=(I-Y_i)^{-1}$ and $a_i=\|w_i\|_2$.
Then $Q_i^\top DQ_i[Z]=2R_i^\top ZR_i$. Using
$R_i^\top[v]_\times R_i=\det(R_i)[R_i^{-1}v]_\times$,
$\det(R_i)=(1+a_i^2)^{-1}$, and
$\|[v]_\times\|_F=\sqrt2\|v\|_2$, its singular values are
$2\sqrt2/(1+a_i^2)$ and $2\sqrt2/\sqrt{1+a_i^2}$, the latter with
multiplicity two. The product Jacobian is block diagonal and
$a_i^2\leq r^2$, proving the rotation-product bound in Table~\ref{tab:additional-radius-bounds}.
For inverse stereographic projection,
$DG(b_i)^\top DG(b_i)=4(1+\|b_i\|_2^2)^{-2}I$.
Each block therefore has singular values $2/(1+\|b_i\|_2^2)$;
using $\|b_i\|_2^2\leq r^2$ proves the oblique bound in Table~\ref{tab:additional-radius-bounds}.
\end{proof}

\begin{proof}[Proof of the complex and generalized Stiefel bounds in Table~\ref{tab:additional-radius-bounds}]
Replace transposes by adjoints in the proof of
the real Stiefel argument above. The same block estimate gives
$\|R^*ZR\|_F\leq(1+b^2)\|P_{\mathfrak m}(R^*ZR)\|_F$,
where $b$ is the spectral norm of the lower-left Cayley block.
The resolvent estimates and quotient differential apply over
$\mathbb C$ with real inner products, proving
the complex Stiefel bound in Table~\ref{tab:additional-radius-bounds}. Whitening is an isometry for $g_B$,
which gives the whitened generalized Stiefel bound in Table~\ref{tab:additional-radius-bounds}.
For a tangent vector $V$ on the ordinary Stiefel manifold,
$\|V\|_F\leq\|V\|_{g_{\mathrm{St}}}\leq\sqrt2\|V\|_F$, so
\[
 \frac{\|V\|_{g_{\mathrm{St}}}}{\sqrt{2\lambda_{\max}(B)}}
 \leq\|B^{-1/2}V\|_F
 \leq\frac{\|V\|_{g_{\mathrm{St}}}}{\sqrt{\lambda_{\min}(B)}}.
\]
The distortion is at most $K_B$, proving
the ambient generalized Stiefel bound in Table~\ref{tab:additional-radius-bounds} by composition. Solving the bound
gives the stated radius.
\end{proof}

\begin{proof}[Proof of \rev{Proposition}~\ref{thm:grassmann-graph}]
Isometric left multiplication reduces to $Q_0=I$.
Set $M=I_p+B^\top B$, $D=I_{N-p}+BB^\top$, and $b=\|B\|_2$.
Orthonormal frames for the subspace and its complement are
\[
 X=\begin{bmatrix}I_p\\B\end{bmatrix}M^{-1/2},
 \qquad X_\perp=\begin{bmatrix}-B^\top\\I_{N-p}\end{bmatrix}D^{-1/2}.
\]
For a parameter variation $C$, horizontal projection eliminates the
normalization derivative and gives
$X_\perp^\top DX(B)[C]=D^{-1/2}CM^{-1/2}$.
Thus $\|DG(B)[C]\|_g=\|D^{-1/2}CM^{-1/2}\|_F$.
Both multipliers have singular values in $[(1+b^2)^{-1/2},1]$, so
$\|C\|_F/(1+b^2)\leq\|DG(B)[C]\|_g\leq\|C\|_F$.
This proves \eqref{eq:grassmann-graph-bounds}; solving the condition-number
bound gives \eqref{eq:grassmann-graph-radius}. Any orthonormal QR frame
differs from $X$ by a right orthogonal factor, whose derivative contributes
only a vertical term. Its horizontal differential has the same norm,
so the bound also applies to the sign-corrected QR implementation.
\end{proof}

All conditioning bounds use the stated coordinate norms and target metrics.
For real Stiefel coordinates, the normal homogeneous metric is
$\|\Delta\|_{g,X}^2=2\|\Delta\|_F^2-\|X^\top\Delta\|_F^2$, and
independent-entry coordinates satisfy $\|Y\|_F^2=2\|y\|_2^2$.
Thus $\|y\|_2\le r_y$ gives the bound $(1+2r_y^2)^2$, multiplied by
$K_B$ for the ambient generalized metric.
For complex Stiefel coordinates,
$\|Y\|_F^2=\|y_{\rm diag}\|_2^2+2\|y_{\rm off}\|_2^2
+2\|y_{\rm rect}\|_2^2$; writing $Y=Sy$ gives
$\kappa_2(J_y)\le\kappa_2(S)(1+\|S\|_2^2r_y^2)^2$ in the normal metric,
with $\|S\|_2,\kappa_2(S)\le\sqrt2$.





\section{Auxiliary stability argument for full-memory BFGS}\label{app:restart-stability}
We supply the extension and matrix estimates used in
\rev{Lemma}~\ref{thm:restart-nonescape}. All constants below are uniform over
centers in a sufficiently small neighborhood of $x^\star$.

\begin{lemma}[Strongly convex extensions and uniform BFGS bounds]\label{lem:restart-extension}
Under \rev{Assumptions}~\ref{ass:restart-global} \rev{and}~\ref{ass:restart-local},
there exist $r>0$ and globally defined $C^{2,1}$ functions
$\widetilde h_x:\R^d\to\R$ agreeing with $h_x$ on $B_r(y_x^\star)$,
with uniform positive Hessian bounds and a uniform Hessian Lipschitz
constant. Standard full-memory BFGS applied to these extensions, from
$y_0=0$ and $H_0=\gamma_xI$ with
$\underline\gamma\le\gamma_x\le\bar\gamma$, using Armijo--Wolfe steps,
satisfies
\begin{equation}\label{eq:extension-uniform-bounds}
 \norm{\widetilde y_{x,j}-y_x^\star}_2
 \le C_Eq^{j/2}\norm{y_x^\star}_2,
 \qquad 0\prec\widetilde H_{x,j}\preceq M_HI,
\end{equation}
where $C_E,M_H<\infty$ and $0<q<1$ are independent of $x$, $j$, and
the particular accepted Wolfe steps. These estimates also hold for every
finite prefix that can be continued as such an auxiliary BFGS sequence.
\end{lemma}
\begin{proof}
Write $u=y-y_x^\star$, $Q_x=\nabla^2h_x(y_x^\star)$, and
$q_x(y)=h_x(y_x^\star)+\tfrac12u^\top Q_xu$.
Choose a fixed smooth cutoff $\chi$ equal to one on $\overline B_1(0)$
and zero outside $B_2(0)$, and define
\[
 \widetilde h_x(y)=q_x(y)+\chi(u/r)\bigl(h_x(y)-q_x(y)\bigr),
\]
with the correction extended by zero outside its support.
Uniform Hessian Lipschitz continuity gives Taylor remainder estimates
of orders $\norm u_2^3$, $\norm u_2^2$, and $\norm u_2$ for the
remainder, gradient, and Hessian. Hence the correction has Hessian norm
at most $C_\chi r$ and a Hessian Lipschitz constant bounded independently
of $r$ and $x$. Taking $r$ small gives uniform constants
$0<a\le b<\infty$ and $K_E<\infty$ such that
\[
 aI\preceq\nabla^2\widetilde h_x(y)\preceq bI,
 \qquad
 \norm{\nabla^2\widetilde h_x(v)-\nabla^2\widetilde h_x(w)}_2
 \le K_E\norm{v-w}_2
\]
for all $y,v,w$. In particular, $y_x^\star$ is the unique minimizer.

The Armijo--Wolfe BFGS estimate of
\cite[Theorem 4.1]{JinJiangMokhtari2024} gives, for $j\ge1$,
\[
 \frac{\widetilde h_x(\widetilde y_{x,j})-\widetilde h_x(y_x^\star)}
 {\widetilde h_x(0)-\widetilde h_x(y_x^\star)}
 \le\left(1-e^{-P_x/j}\frac{2c_1(1-c_2)a}{b}\right)^j,
 \quad P_x=\Psi(B_0/b),
\]
where $B_0=H_0^{-1}$ and
$\Psi(A)=\operatorname{tr}(A)-\log\det A-d$.
(The zero initial-gap case is already stationary.)
The bounds on $\gamma_x$ give $P_x\le\bar P<\infty$.
Thus $q=1-e^{-\bar P}2c_1(1-c_2)a/b\in(0,1)$ and strong convexity yield
the first estimate in \eqref{eq:extension-uniform-bounds}, with
$C_E=\sqrt{b/a}$. In particular,
\begin{equation}\label{eq:extension-summable}
 \sum_{j=0}^\infty\norm{\widetilde y_{x,j}-y_x^\star}_2
 \le\frac{C_E}{1-\sqrt q}\norm{y_x^\star}_2.
\end{equation}

For the matrix bound, use the normalized variables
\[
 A_j=Q_x^{-1/2}B_jQ_x^{-1/2},\qquad
 \widehat s_j=Q_x^{1/2}s_j,\qquad
 \widehat z_j=Q_x^{-1/2}z_j,
\]
where $z_j$ is the exact gradient difference for $\widetilde h_x$.
The BFGS trace and determinant identities give
\begin{align*}
 \Psi(A_{j+1})-\Psi(A_j)
 &=\frac{\norm{\widehat z_j}_2^2}{\widehat s_j^\top\widehat z_j}
   -\frac{\norm{A_j\widehat s_j}_2^2}{\widehat s_j^\top A_j\widehat s_j}
   -\log\frac{\widehat s_j^\top\widehat z_j}
                    {\widehat s_j^\top A_j\widehat s_j}\\
 &\le\frac{\norm{\widehat z_j}_2^2}{\widehat s_j^\top\widehat z_j}
       -1-\log\frac{\widehat s_j^\top\widehat z_j}
                         {\norm{\widehat s_j}_2^2}.
\end{align*}
The inequality uses Cauchy--Schwarz and $\log t-t\le-1$.
The global Hessian bounds keep
$\widehat s_j^\top\widehat z_j/\norm{\widehat s_j}_2^2$
uniformly positive and $\norm{\widehat z_j}_2/\norm{\widehat s_j}_2$
uniformly bounded. Hessian Lipschitz continuity also gives
\[
 \frac{\norm{\widehat z_j-\widehat s_j}_2}
      {\norm{\widehat s_j}_2}
 \le C_K\max\{\norm{\widetilde y_{x,j}-y_x^\star}_2,
                 \norm{\widetilde y_{x,j+1}-y_x^\star}_2\}.
\]
The preceding scalar expression vanishes when $\widehat z_j=\widehat s_j$
and is uniformly Lipschitz on these bounded, positive-denominator ranges.
Consequently,
\[
 \Psi(A_{j+1})-\Psi(A_j)
 \le C_P\bigl(\norm{\widetilde y_{x,j}-y_x^\star}_2+
                  \norm{\widetilde y_{x,j+1}-y_x^\star}_2\bigr).
\]
Together with \eqref{eq:extension-summable} and the uniform bound on
$\Psi(A_0)$, this bounds $\Psi(A_j)$ uniformly. Since
$\Psi(A_j)=\sum_i(\lambda_i(A_j)-\log\lambda_i(A_j)-1)$,
its sublevel sets have eigenvalues uniformly bounded away from zero and
infinity. Thus $H_j=B_j^{-1}\preceq M_HI$ uniformly.
A finite admissible prefix can always be continued for the globally
strongly convex extension by ordinary Wolfe steps, so the same constants
apply to every such prefix, independently of how it is continued.
\end{proof}

\end{document}